\documentclass[11pt]{amsart}
\usepackage{latexsym,  amsfonts, amssymb, graphicx, a4, epic, eepic, epsfig,
setspace, lscape}

\usepackage{color}

\def\inst#1{$^{#1}$}

\newtheorem{theorem}{Theorem}[section]
\newtheorem{lemma}[theorem]{Lemma}
\newtheorem{proposition}[theorem]{Proposition}

\newtheorem{definition}[theorem]{Definition}
\newtheorem{corollary}[theorem]{Corollary}
\newtheorem{remark}[theorem]{Remark}

\def \P {{\mathbb P}}
\def \E {{\mathbb E}}

\def \L {{\Lambda}}

\def \s {{\sigma}}

\def \g {{\gamma}}
\def \t {{\tau}}

\def \d {{\delta}}
\def \p {{\pi}}
\def \x {{\xi}}

\def\ra{\rightarrow}
\def\Proof{{\sl Proof.}\quad}

\newcommand{\be}[1]{\begin{equation}\label{#1}}
\newcommand{\ee}{\end{equation}}

\newcommand{\bl}[1]{\begin{lemma}\label{#1}}
\newcommand{\el}{\end{lemma}}

\newcommand{\br}[1]{\begin{remark}\label{#1}}
\newcommand{\er}{\end{remark}}

\newcommand{\bt}[1]{\begin{theorem}\label{#1}}
\newcommand{\et}{\end{theorem}}

\newcommand{\bd}[1]{\begin{definition}\label{#1}}
\newcommand{\ed}{\end{definition}}

\newcommand{\bcl}[1]{\begin{claim}\label{#1}}
\newcommand{\ecl}{\end{claim}}

\newcommand{\bp}[1]{\begin{proposition}\label{#1}}
\newcommand{\ep}{\end{proposition}}

\newcommand{\bc}[1]{\begin{corollary}\label{#1}}
\newcommand{\ec}{\end{corollary}}

\newcommand{\bi}{\begin{itemize}}
\newcommand{\ei}{\end{itemize}}

\newcommand{\ben}{\begin{enumerate}}
\newcommand{\een}{\end{enumerate}}

\def \qed {${\square\hfill}$\newline}
\def \Z {{\mathbb Z}}

\def \N {{\mathbb N}}
\def \P {{\mathbb P}}
\def \E {{\mathbb E}}

\begin{document}

\title[Fast mixing for the 2d Ising model]{Fast mixing for the low temperature 2d Ising model through irreversible parallel  dynamics}

\author{
Paolo Dai Pra\inst{1}\and
Benedetto Scoppola\inst{2} \and
Elisabetta Scoppola\inst{3}}


\maketitle
\begin{center}
{\footnotesize
\vspace{0.3cm} \inst{1}  Dipartimento di Matematica,
University of Padova\\
Via Trieste 63, 35121 Padova, Italy\\
\texttt{daipra@math.unipd.it}\\

\vspace{0.3cm}\inst{2}  Dipartimento di Matematica, University of Roma
``Tor Vergata''\\
Via della Ricerca Scientifica - 00133 Roma, Italy\\
\texttt{scoppola@mat.uniroma2.it}\\

\vspace{0.3cm} \inst{3} Dipartimento di Matematica e Fisica, University of Roma
``Roma Tre''\\
Largo San Murialdo, 1 - 00146 Roma, Italy\\
\texttt{scoppola@mat.uniroma3.it}\\ }

\end{center}


\begin{abstract}
We study metastability and mixing time for a non-reversible probabilistic cellular automaton. With a suitable choice of the parameters, we first show
that the stationary distribution is close in total variation to a low temperature Ising model. Then we prove that both the mixing time and the time to exit a 
metastable state grow polynomially in the size of the system, while this growth is exponential in reversible dynamics. In this model, non-reversibility, parallel updatings and a suitable choice of boundary conditions combine to produce an efficient dynamical stability.
\end{abstract}


\tableofcontents

\section{Introduction}
\label{Intro}

In this paper we consider a discrete-time stochastic dynamics for a spin system at low temperature, in which high mobility of parallel
updating and asymmetry of the interaction combine to produce efficient dynamical stability and fast convergence to equilibrium. 

%

The control of the convergence to equilibrium of irreducible Markov Chains (MC)  is particularly interesting when the invariant measure is strongly polarized,
for instance in MC describing large scale ferromagnetic systems at low temperature.
Indeed  in the region of parameters where the system  exhibits coexistence of more phases,
the problem of the control of the convergence to equilibrium of the MC describing the system becomes
strictly related to the problem of {\em metastability}, since the tunneling between different phases is necessary
to reach equilibrium. This tunneling time usually is exponentially divergent in the size of the problem so that the
convergence to equilibrium in these cases is exponentially slow. See \cite{fabio} for a beautiful review on 
this problem.

We briefly recall the well known Ising model in 2-d  in order to explain in more detail the problem.

Let $L$ be a positive integer, and $\L := \left(\Z/L\Z \right)^2$ be the two dimensional discrete torus. 
Consider the standard Ising model on $\L$ without external field with {\em spin configurations} $\s = (\s_x)_{x \in \L} \in \mathcal{S}:= \{-1,1\}^{\L}$ and with Hamiltonian
\be{HG}
H(\s)=-\sum_{(x,y)}J\s_x\s_y
\ee
where $J>0$ and the sum is on neighboring sites in $\L$.   Denote by
$\p_{G}$ its Gibbs measure
\be{Gm}
\p_{G}(\s)=\frac{e^{-H(\s)}}{Z_{G}},\qquad Z_{G}=\sum_{\s\in \mathcal{S}}e^{-H(\s)}.
\ee

A popular discrete-time MC, reversible w.r.t. this Gibbs measure, is given by the following algorithm: at each time $t$
a point $x \in \L$ is chosen with uniform probability; all spins $\s_y$, $y \neq x$ are left unchanged, while $\s_x$ is {\em flipped} with 
probability
\[
\exp\left[ -(H(\s^x) - H(\s))^+ \right],
\]
where $\s^x$ is the configuration obtained by $\s$ by flipping $\s_x$ and, for a real number $a$, 
$a^+ := \max(a,0)$.
Denote by $\P_{\s}^t$ the probability distribution of the process at time $t$
starting from $\s$ at time $0$, and by $\mathcal{P}$ the transition matrix $(\P_{\s}^1(\eta))_{\s,\eta \in \mathcal{S}}$.

Different quantities can be used to control the convergence to equilibrium of MC's;
the most popular is the {\it mixing time}
\be{mixtime}
T_{mix}:=\min\Big\{t>0;\; d(t)\le \frac{1}{e}\Big\}
\ee
where $d(t)$ is the maximal distance in total variation between the distribution at time $t$
and the invariant measure
$$d(t)=\sup_\s\Vert \P_{\s}^t-\p_G\Vert_{TV}.$$

%

For the Glauber dynamics defined above when  the interaction constant $J$ is so large
  that the Gibbs measure $\p_{G}$ is  nearly concentrated on the
configurations ${\bf +1}$  and ${\bf -1}$, with all spins +1 and all spins -1 respectively,
it is possible to prove  that $T_{mix}$ diverges exponentially in $L$.
 This result is due to the presence of a rather tight
``bottleneck" in the state space. Indeed starting for instance from ${\bf -1}$, in order to relax to
equilibrium the dynamics has to reach a neighborhood of the opposite minimum ${\bf +1}$, crossing
the set of configurations with zero magnetization which has a small Gibbs measure.
 In other words the system is trapped for a very long time near the
configuration ${\bf -1}$, and only after many attempts to leave this trap,  a +1 droplet is nucleated and grows up to
reach the bottleneck, i.e.,  the set of configurations of zero magnetization. This mechanism
is typical in metastability and produces a large relaxation time.

If the relaxation time is exponentially large, the MC given by the Glauber dynamics is
not an efficient way for sampling from the Gibbs measure $\p_G$ for large systems.
A possible way to bypass this problem is the following: for each size $L$ of the system, 
we construct a MC whose invariant measure $\pi$ is  close to $\pi_G$, in the sense that 
$\Vert \p-\p_{G}\Vert_{TV}$ converges to zero as  $L \ra +\infty$, and such that its mixing time 
grows polynomially in $L$. We call this a {\em asymptotically polynomial approximation scheme}; this notion is weaker but closely related to that
of {\em fully polynomial randomized 
approximation scheme} (FPRAS) introduced in theoretical computer science (see [Jerrum, Sinclair]).


In this paper we present two independent results;  their combination provide
an asymptotically polynomial approximation scheme for the 2d Ising model.

More precisely we introduce a modification of the above MC in which:
\begin{itemize}
\item
all spins are simultaneously updated;
\item
the updating of the spin $\s_x$ only depends on the spin values at the previous time of its
South and West nearest neighbors; this makes the dynamics non reversible.
\end{itemize}
This dynamics is a 
 probabilistic cellular automaton (PCA) for which the invariant measure $\p_{PCA}$ can be found
without a  detailed balance condition, but using the notion of {\it weak balance condition}
discussed in \cite{LS} also known in the literature as {\it dynamical reversibility} (see for instance
\cite{A}). By using the ideas developed in \cite{DSS} we can control
the total variation distance between the Gibbs measure $\p_{G}$ and 
$\p_{PCA}$. This is the content of Theorem \ref{theo}.

In the second theorem we study the convergence to equilibrium of this parallel dynamics. 
The key step is an estimate on the tunneling time between ${\bf -1}$ and ${\bf 1}$. This estimate is obtained 
by using some of the basic ideas developed in the context of metastability. The main point is concerned 
with the separation of time scales. The general idea is the following: the energy landscape determines a sequence
$\mathcal{S} = \mathcal{S}_0 \supset \mathcal{S}_1 \supset \cdots \supset \mathcal{S}_n = \{{\bf -1}, {\bf 1} \}$
of nested subsets of $\mathcal{S}$ is such a way that for $k \geq 1$ a time scale $T_k$ is associated to each $\mathcal{S}_k$
in the following sense: the dynamics need a time of order $T_k$ to leave $\mathcal{S}_k$, but a much smaller time to return in $\mathcal{S}_k$ after having
left it; moreover, $T_k$ is much smaller than $T_{k+1}$. This allows to define an effective {\em renormalized} dynamics on $\mathcal{S}_k$ which evolves 
at time-scale $T_k$, and which consists of the successive returns in $\mathcal{S}_k$. See for instance  \cite{S},\cite{OV} and \cite{L} for more details
on such a renormalisation procedure. Iterating this strategy on larger and larger time scales $t_0<t_1<... <t_n$ 
one arrives to the situation in which $\mathcal{S}_n$ is given just by the absolute minima of the energy.
In this case the corresponding renormalized process is a very elementary two states process with a tunneling time
$\t({\bf -1},{\bf +1})$ given by an exponential random variable with mean given by the inverse of the
transition probability $({\bf -1},{\bf +1})$ of the renormalized chain on  $\mathcal{S}_n$.

We do not  completely develop this analysis for our PCA dynamics, but we will use the main ideas of 
separation of time scales and corresponding 
reduction of the state space in order to control the mean tunneling time, and, with this, the mixing time
of the PCA.
Exploiting  the complete asymmetry of the interaction (only SW), the simultaneous updating and the periodic boundary conditions, we observe that configurations with the same spin
on a NW-SE diagonal are stable on the time scale of order 1, just moving in the NE direction.
Playing on the difference of time scales involved in the process, we can tune the parameters of the dynamics
in order to describe the evolution between diagonal configurations in terms of a 1d nearly symmetric Random Walk,
producing a tunneling time which is polynomial in the size of $\L$.

In Section \ref{sec:model} we define the model in details, and state our main results. Section \ref{sec:lowT} is devoted to the analysis of the invariant 
measure of the PCA, and its relations with the Ising model. Some fundamental facts on time scale separation for the PCA is presented in Section \ref{sec:low temp}, while Section \ref{sec:mix} contains the key estimate on the tunneling time.

\section{The model and the results} \label{sec:model}

\subsection{The model}

On the same space of configurations $\mathcal{S}:= \{-1,1\}$ discussed in the Introduction
for the Ising model we want to construct a
Markov chain given in terms of a completely asymmetric interaction as follows.
For $x = (i,j) \in \L = \left(\Z/L\Z \right)^2$, we introduce the following notation for its nearest neighbors:
\be{xrd}
x^u := (i,j+1) \ \ \ x^r := (i+1,j) \ \ \ x^d := (i, j-1) \ \ \ x^l := (i-1,j)
\ee
where sums and difference has to be meant mod. $L$. Given a { spin configuration} $\s = (\s_x)_{x \in \L} \in \mathcal{S}$, for typographical reasons we write $\s_x^u$ for $\s_{x^u}$, and similarly for the other nearest neighbors of $x$. Consider the  discrete-time Markov chain on $\mathcal{S}$, whose transition matrix is given by
\be{Pst}
P({\s,\t}) := \frac{e^{-H(\s,\t)}}{\sum_{\s' \in \mathcal{S}} e^{-H(\s,\s')}},
\ee
where $H(\s,\t)$ is the following asymmetric Hamiltonian, defined on pairs of configurations:
\be{Hst}
\begin{split}
H(\s,\t)  & := - \sum_{x \in \L} \left[ J \s_x(\t_x^u + \t_x^r) + q \s_x \t_x \right] \\ & = - \sum_{x \in \L} \left[ J \t_x(\s_x^d + \s_x^l) + q \s_x \t_x \right]
\end{split}
\ee
and $J,q>0$ are given parameters. In what follows we set 
\be{Zs}
Z_{\s} := \sum_{\s' \in \mathcal{S}} e^{-H(\s,\s')}.
\ee

Some basic facts on this Markov chain are grouped in the next Proposition (see \cite{LS} for more details) motivating the name {\it Probabilistic Cellular Automata} (PCA) for this dynamics.
\bp{prop:basic}
\ben
\item
$P({\s,\tau})$ is of the following product form:
\[
P({\s,\tau}) = \prod_{x \in \L} p_x(\tau_x|\s)
\]
where
\[
p_x(\tau_x|\s) := \frac{\exp\left\{ \tau_x\left[J(\s_x^d + \s_x^l) + q \s_x\right]\right\}}{2 \cosh(J(\s_x^d + \s_x^l) + q \s_x)}.
\]
\item
$H(\s,\t)\not=H(\t,\s)$ but the following {\em weak symmetry condition} holds
\[
\sum_{\tau \in \mathcal{S}} e^{-H(\s,\tau)} = \sum_{\tau \in \mathcal{S}} e^{-H(\tau,\s)} .
\]
\item
The Markov chain is irreversible with a unique stationary distribution $\pi_{PCA}$ given by
\[
\pi_{PCA}(\s) := \frac{Z_{\s}}{Z_{PCA}},
\]
with $Z_{PCA} := \sum_{\s} Z_{\s}$.
\een
\ep
\Proof
The statement in (1) amounts to a straightforward computation; in particular, it implies irreducibility of the chain, which therefore has a unique stationary distribution. The statement in (3) thus follows readily from (2), that is the only nontrivial point to show. Note that
\be{1}
\begin{split}
\sum_{\tau \in \mathcal{S}}  e^{-H(\s,\tau)} & =   2^{|\L|} \prod_{x \in \L}  \cosh(J(\s_x^d + \s_x^l) + q \s_x)) \\
 \sum_{\tau \in \mathcal{S}} e^{-H(\tau,\s)} & =  2^{|\L|} \prod_{x \in \L}  \cosh(J(\s_x^u + \s_x^r) + q \s_x)).
 \end{split}
 \ee
Denote by $\L^* := \{ \{x,y\}: \,  \x,y \in \L, \, |x-y| = 1\}$ the set of {\em bonds} in $\L$. Note that $|\L^*| = 2L^2$. For $\s \in \mathcal{S}$, we let 
\be{peierls}
\g(\s) := \{ \{x,y\} \in \L^* : \, \s_x \neq \s_y \}
\ee
be the {\em Peierls contour} associated to $\s$. The following identities are immediately checked:
\[
\cosh(J(\s_x^d + \s_x^l) + q \s_x)) = \left\{
\begin{array}{ll}
\cosh(2J+q) & \mbox{if } \{x,x^d\} \not\in \g(\s), \, \{x,x^l\} \not\in \g(\s) \\
\cosh(2J-q) & \mbox{if } \{x,x^d\} \in \g(\s), \, \{x,x^l\} \in \g(\s) \\
\cosh(q) & \mbox{otherwise}.
\end{array}
\right.
\]
So, if we let
\[
n_{dl} = n_{dl}(\s) := \left| \{ x \in \L: \, \{x,x^d\} \in \g(\s), \, \{x,x^l\} \in \g(\s) \} \right|,
\]
using \eqref{1} we obtain
\be{2}
\sum_{\tau \in \mathcal{S}}  e^{-H(\s,\tau)}  = 2^{L^2} [\cosh(2J-q)]^{n_{dl}}[\cosh(q)]^{|\g(\s)| - 2 n_{dl}}[\cosh(2J+q)]^{L^2 -|\g(\s)| + n_{dl}}.
\ee
With the same argument, defining
\be{3}
n_{ur} =n_{ur}(\s) := \left| \{ x \in \L: \, \{x,x^u\} \in \g(\s), \, \{x,x^r\} \in \g(\s) \} \right|,
\ee
we obtain
\be{4}
\sum_{\tau \in \mathcal{S}}  e^{-H(\tau,\s)}  = 2^{L^2} [\cosh(2J-q)]^{n_{ur}}[\cosh(q)]^{|\g(\s)| - 2 n_{ur}}[\cosh(2J+q)]^{L^2 -|\g(\s)| + n_{ur}}.
\ee
The conclusion now follows from the observation that, for every $\s \in \mathcal{S}$, the identity $n_{dl}(\s) = n_{ur}(\s)$ holds. This can be shown, for instance, by induction on $n^+(\s)$, where $n^+(\s)$ denotes the number of spins equal to $+1$ in $\s$. If $n^+(\s) = 0$ the statement is obvious. For $n^+(\s) = n >0$, let $x \in \L$ be such that $\s_x = +1$, and let $\s^x$ the configuration obtained from $\s$ by flipping the spin at $x$. By considering all possible spin configuration in the $3 \times 3$ square centered at $x$, one checks that $n_{dl}(\s^x) - n_{ur}(\s^x) = n_{dl}(\s) - n_{ur}(\s)$. Since $n^+(\s^x) = n^+(\s) -1$, the proof is completed.
\qed


\subsection{The results}

We are interested in the limit $L\to\infty$ and in the low temperature ($J$ large) regime  defined as follows. 
\bd{lt-regime}
	The {\em low temperature regime with parameters  $k$ and $c$} corresponds to the
following choice
\be{regime}
 J = J(L) = k \log L\qquad q = q(L) = c \frac{\log L}{L}
 \ee
  \ed
We state here our two main results. The first concerns the relation between the two considered models,
controlling the distance in total variation between the Gibbs measure of the
symmetric standard Ising model and the stationary distribution of the asymmetric PCA. The numerical constants
appearing in the statements of the theorems are not optimized.

\begin{theorem}\label{theo}
In the low temperature regime with parameter $k$ and $c$,  there is a constant $C>0$ such that
\be{tvdist}
\Vert\pi_{PCA}-\pi_G\Vert_{TV}\le C \left(\frac{1}{L^{\frac{c}{2}-1}} + \frac{1}{L^{2k-2}} \right).
\ee
\end{theorem}

The second result is the control of the convergence to equilibrium of the PCA proving that the mixing time of the parallel dynamics is polynomial in $L$.

\begin{theorem}\label{mix}
In the low temperature regime with parameter $k$ and $c$ such that $c > \frac12$ and $k-4c >4$, we have
$$
\lim_{L\to\infty} d_{PCA}(L^{8k})=0
$$
where 
$$
d_{PCA}(t)=\sup_\s||P^t(\s,.)-\p_{PCA}(.)||_{TV}
$$
\end{theorem}

Theorems \ref{theo} and \ref{mix} imply that the Markov chain defined in (\ref{Pst})
provides a asymptotically polynomial approximation scheme for the  Ising model on the 2d torus.

\br{rbc}
As mentioned in the introduction, the periodic boundary conditions play a crucial role in the proof of 
Theorem  \ref{theo}.
\er

\br{rth} 
There is another example, see \cite{fabio2}, 
of rapid mixing of a Markov chain having 
as stationary measure the Gibbs measure of the low temperature Ising model. This example is the
Swendsen-Wang dynamics. As in our case, such dynamics is fast because it allows the possibility to update 
in a single step of the Markov chain a large amount of spins.
However, as far as we know, this is the first case in literature of a fast {\it irreversible} dynamics
based on the idea of the PCA. In
particular it seems that  the ingredient of the irreversibility combined with parallelism  is quite crucial in order to obtain 
the fast mixing. 
\er

\section{The relation between Ising Gibbs measure and PCA stationary measure
at low temperature} \label{sec:lowT}

We prove in this section Theorem \ref{theo}.

We use the representation introduced in \cite{DSS}.
Note first of all that 
\be{ww}
\begin{split}
Z_\s & =\sum_\t e^{-\sum_{x} [J(\sigma_x^d+\sigma_x^l)+q\s_x]\t_x} \\ & =
e^{q|\L|}\sum_{I\subset \L}e^{\sum_{(x,y)}J\s_x\s_y-2\sum_{x\in I}J(\s_x\s_x^u+\s_x\s_x^r)-2q|I|} \\
&  = e^{q|\L|}w^{G}(\s)\prod_{x\in \L}(1+\delta\phi_x)
\end{split}
\ee
where we have used $\d=e^{-2q}$, 
\[
w^{G}(\s) = e^{-H(\s)},
\]
 and 
$$\phi_x=e^{-2J(\s_x\s_x^u+\s_x\s_x^r)}.$$
We will call
\be{f}
f(\sigma)=\prod_{x\in \L}(1+\delta\phi_x).
\ee
It easily follows that
\be{pgppca}
\p^{PCA}(\s)=\p^{G}(\s)\frac{f}{\p^G(f)}
\ee

\noindent
We have then
\begin{equation}\label{2.1}
\Vert\pi_{PCA}-\pi_G\Vert_{TV}=\pi_G\left[\left|\frac{f}{\pi_G(f)}-1\right|\right]
\end{equation}

\noindent
Write now  the Gibbs measure in terms of Peierls contours (see \eqref{peierls}):
$$\pi_G(\s)=\frac{e^{-2Jl(\s)}}{Z_G}
$$
where $l(\s) := |\g(\s)|$ is the total length of the Peierls contours of the configuration $\s$.

\noindent
Let ${\bf 1}$ be the configuration with $\s_x=1$ for all $x$.

\noindent
Normalizing the value of $f(\s)$ with the value $f({\bf 1})$, 
which is a constant ineffective in the evaluation of (\ref{2.1}),
the expression of $f(\s)$ can be written as (see also \eqref{3})
\be{eff}
f(\s)=\left[\frac{(1+\d e^{4J})}{(1+\d e^{-4J})}\right]^{n_{ur}(\s)}
\left[\frac{(1+\d)}{(1+\d e^{-4J})}\right]^{l(\s)-2n_{ur}(\s)}.
\ee
where we have simply observed that 
\[
\s_x\s_x^u+\s_x\s_x^r = \left\{ \begin{array}{ll} 2 & \mbox{if } (x,x^u) \in \g(\s), (x,x^r) \in \g(\s) \\ 0 & \mbox{if } (x,x^u) \not\in \g(\s), (x,x^r) \not\in \g(\s) \\ 1 & \mbox{otherwise.} \end{array} \right.
\]

\noindent
Note that with this normalization 
$f({\bf 1})=1$ obviously holds.

\noindent
Let us first give an upper bound of $\pi_G(f)$.
We can write
\[
\begin{split}
\pi_G(f) & =\frac{1}{Z_G}\sum_\s  \left[e^{-4J}\frac{(1+\d e^{4J})}{(1+\d e^{-4J})}\right]^{n_{ur}(\s)}
\left[e^{-2J}\frac{(1+\d)}{(1+\d e^{-4J})}\right]^{l(\s)-2n_{ur}(\s)} \\ 
&
\le \frac{1}{Z_G}\sum_\s\left[\d+e^{-4J}\right]^{n_{ur}(\s)}\left[2e^{-2J}\right]^{l(\s)-2n_{ur}(\s)}. 
\end{split}
\]
To give estimates for this last sum, we use again Peierls contours. We say that a pair of adjacent bonds $(x,x^u), (x,x^r)$ both belonging to $\g(\s)$ form a {\em ur-elbow}. Note that the only closed paths in $\L^*$ exclusively consisting of ur-elbows is necessarily union of complete diagonals (actually of a even number of diagonals, for the contour to correspond to a spin configuration). Any contour $\g = \g(\s)$ can be decomposed as $\g = \g_D \cup \g_{ND}$, where $\g_D$ only contains complete diagonals, while $\g_{ND}$ has no complete diagonal. Observe that $l(\s)-2n_{ur}(\s) = 0 \ \iff \ \g_{ND}(\s) = \emptyset$. Now, for any fixed $m \geq 0$ we obtain an upper bound for the contribution of all configurations $\s$ such that $l(\s)-2n_{ur}(\s) =  m$. We can write
\[
\begin{split}
A(m) & := \sum_{\s: l(\s)-2n_{ur}(\s) =  m} \left[\d+e^{-4J}\right]^{n_{ur}(\s)}\left[2e^{-2J}\right]^{l(\s)-2n_{ur}(\s)}  \\ & =2 \sum_{\g : |\g|-2n_{ur}(\g) =  m}
\left[\d+e^{-4J}\right]^{n_{ur}(\g)}\left[2e^{-2J}\right]^{m},
\end{split}
\]
where the factor $2$ come from the fact that there are exactly two configurations for each contour.
Observe now that $e^{-2J} = 1/L^{2k}$ while
\[
\d+e^{-4J} = e^{-2c \frac{\log L}{L}} + e^{-4k \log L} \leq 1-c \frac{\log L}{L} + \frac{1}{L^{4k}}  \leq 1- \frac{c}{2} \frac{\log L}{L} < 1,
\]
for $L$ sufficiently large. Thus, using the decomposition $\g = \g_D \cup \g_{ND}$,
\[
\begin{split}
A(m) & \leq  2\sum_{\g : |\g|-2n_{ur}(\g) =  m} \left(1- \frac{c}{2} \frac{\log L}{L} \right)^{n_{ur}(\g_D)} \left( \frac{2}{L^{2k}} \right)^m \\
& \leq  2\left( \frac{2}{L^{2k}} \right)^m N_m \sum_{\g: \g_{ND} = \emptyset} \left(1- \frac{c}{2} \frac{\log L}{L} \right)^{\frac{|\g|}{2}},
\end{split}
\]
where 
\[
N_m := \left| \{ \g : \g_D = \emptyset, \, |\g|-2n_{ur}(\g) =  m\} \right|.
\]
A very rough upper bound for $N_m$ can be obtained as follows. We first place the $m$ bonds not belonging to a ur-elbow (we have at most $(2L^2)^m$ different choices); call $\tilde{\g}_{ND}$ the resulting set of bonds. We then place an arbitrary number of ur-elbows, with the constraint that the endpoints of a connected sequence of NE elbows must coincide with two of the $2m$ endpoints of $\tilde{\g}_{ND}$. Moreover, for any endpoint $x$ of $\tilde{\g}_{ND}$ there are at most two connected sequences of ur-elbows which connect $x$ to exactly one endpoint of $\tilde{\g}_{ND}$. Thus, sequences of ur-elbows can be placed in at most $4^{2m}$ different ways. This yields
\[
N_n \leq \left( 32 L^2 \right)^m.
\]
To complete the upper bound for $A(m)$, we need to estimate
\[
\sum_{\g: \g_{ND} = \emptyset} \left(1- \frac{c}{2} \frac{\log L}{L} \right)^{\frac{|\g|}{2}}.
\]
Since such diagonal contours are just union of complete diagonals, and each complete diagonal has length $2L$, for $L$ sufficiently large we have
\[
\begin{split}
\sum_{\g: \g_{ND} = \emptyset} \left(1- \frac{c}{2} \frac{\log L}{L} \right)^{\frac{|\g|}{2}} & \leq \sum_{l \geq 0} \binom{L}{l} \left(1- \frac{c}{2} \frac{\log L}{L} \right)^{lL}  \\
& = \left[1+ \left(1- \frac{c}{2} \frac{\log L}{L} \right)^L \right]^L  \leq 1+\frac{2}{L^{\frac{c}{2}-1}}.
\end{split}
\]
Thus we have
\[
A(m) \leq 2 \left(1- \frac{c}{2} \frac{\log L}{L} \right)^{\frac{|\g|}{2}} \left( \frac{64}{L^{2k-2}} \right)^m .
\]
Summing up, using also the obvious fact that $Z_G = \sum_{\s} e^{-2J l(\s)} > 2$,
we can choose $C>0$ such that for $L$ large enough:
\be{pi(f)}
\begin{split}
\pi_G(f) &  \leq \frac{1}{Z_G}\sum_{m \geq 0} A(m) \leq \left( 1+\frac{2}{L^{\frac{c}{2}-1}} \right) \left[ \sum_{m \geq 0} \left( \frac{64}{L^{2k-2}} \right)^m \right] \\
 & \le 1+ \frac{C}{L^{\frac{c}{2}-1}} + \frac{C}{L^{2k-2}}.
\end{split}
\ee
Comparing \eqref{eff} with \eqref{pi(f)}, and using the fact that $\d \simeq 1$ (in particular $\d \geq \frac12$) for large $L$, one realizes that
 $f(\s) > \pi_G(f)$ for all configurations different from 
$\pm {\bf 1}$. This is evident for $n_{ur}(\s)>0$; for $n_{ur}(\s)=0$ we have that 
if $l(\s)>0$, then
$l(\s)\ge L$, giving $f\ge (1+1/4)^L$.
By this observation
$$\pi_G\left[\left|\frac{f}{\pi_G(f)}-1\right|\right]=\frac{2}{\pi_G(f)}\sum_{\s:f(\s)<\pi_G(f)}
\frac{e^{-H(\s)}}{Z_G}[\pi_G(f)-f(\s)]$$
and the sum contains actually only the two configurations $\s=\pm\bf 1$, such that $f(\s)=1$.

\noindent 
Hence we have
\begin{equation}\label{2.2}
\Vert\pi_{PCA}-\pi_G\Vert_{TV}=\pi_G\left[\left|\frac{f}{\pi_G(f)}-1\right|\right]\le
\frac{2}{\pi_G(f)}[\pi_G(f)-1]=2(1-\frac{1}{\pi_G(f)})\le 2C \left(\frac{1}{L^{\frac{c}{2}-1}} + \frac{1}{L^{2k-2}} \right).
\end{equation}

\noindent
Inserting \eqref{pi(f)} in (\ref{2.2}), and using \eqref{2.1}, we complete the proof of the theorem.

\section{PCA at low temperature} \label{sec:low temp}

\subsection
{Realization through random numbers}

In what follows it will be useful to realize the Markov chain described above using uniformly distributed random numbers. Let $\{U_x(n): \, x \in \L, \, n \geq 1\}$ be a family of i.i.d. random variables, uniformly distributed in $(0,1)$, defined in some probability space $(\Omega,\mathcal{A},P)$. Given the initial configuration $\s(0)$, define recursively $\s(n+1)$ by: $\s_x(n+1) = 1$ if and only if {\em one} of the following conditions holds:
\bi
\item[(A)]
$\s_x^d(n) = \s_x^l(n) = 1$ and $U_x(n+1) \leq \frac{e^{2J + q \s_x(n)}}{2 \cosh(2J + q \s_x(n))}$;
\item[(B)]
$\s_x^d(n) = \s_x^l(n) = -1$ and $U_x(n+1) \leq \frac{e^{-2J + q \s_x(n)}}{2 \cosh(-2J + q \s_x(n))}$;
\item[(C)]
$\s_x^d(n) = -\s_x^l(n) $ and $U_x(n+1) \leq \frac{e^{q \s_x(n)}}{2 \cosh(q \s_x(n))}$,
\ei
while $\s_x(n+1) = -1$ otherwise.

\br{monotonia}
Note that with this construction of the process it is immediate to see that the Markov chain
preserves the componentwise partial order on configuration. Coupling the processes $(\s(n))_{n\in\N}$ starting at
$\s$ and $(\s'(n))_{n\in\N}$ starting at $\s'$ by using the same realization  of uniform variables
 $\{U_x(n): \, x \in \L, \, n \geq 1\}$ we have that if $\s\le\s'$ in the sense that $\s_x\le\s'_x$ for all
 $x\in \L$ then $\s(n)\le\s'(n)$ for each time $n\ge 0$.

\er
\subsection
{Zero-temperature dynamics}

In the low temperature regime considered in this paper, updatings of type (A) or, symmetrically, those for which $\s_x^d(n) = \s_x^l(n) = -1$ $ \mapsto$  $\s_x(n+1) = -1$, are {\em typical}, as they occur with probability $\frac{e^{2J \pm q}}{2 \cosh(2J\pm q} \simeq 1$; conversely updatings of type (B), or those for which  $\s_x^d(n) = \s_x^l(n) = 1$ $ \mapsto$  $\s_x(n+1) = -1$, are {\em atypical}, as they occur with probability $\frac{e^{-2J \pm q}}{2 \cosh(2J \pm q)} \simeq 0$. Finally, updatings of type (C), or those for which  $\s_x^d(n) = -\s_x^l(n) = 1$ $ \mapsto$  $\s_x(n+1) = -1$, are {\em neutral}, as they occur with probability $ \frac{e^{\pm q}}{2 \cosh(q)} \simeq \frac12$. 

Let $N > 0$ be a given (large) time. In next section it will be useful to rule out events of very small probability. For instance, given a time $N>0$, we can ``force'' the system to perform no atypical updating up to time $N$. To this aim, we can define 
\be{no-aty}
S := \min\left\{n \geq 1: \, \exists \, x \mbox{ such that } U_x(n)  \not\in \left( \frac{e^{-2J+q}}{2\cosh(2J-q)}, 1- \frac{e^{-2J+q}}{2\cosh(2J-q)}\right) \right\},
\ee
and condition to the event $\{S>N\}$. Note that under $\P(\cdot | S>N)$ the random numbers $\{U_x(n): \, x \in \L, \, 1 \leq n \leq N\}$ are i.i.d., uniformly distributed on $\left( \frac{e^{-2J+q}}{2\cosh(2J-q)}, 1- \frac{e^{-2J+q}}{2\cosh(2J-q)}\right)$. Thus, $(\s(n))_{n=0}^N$ is a homogeneous Markov chain also under $\P(\cdot | S>N)$, for which only typical and neutral transitions are allowed. This conditioned
dynamics is often called  the  {\it zero-temperature dynamics}
corresponding, for the inverse temperature parameter $J$, to the limit $J\to\infty$. 

Note also that, if $A$ is an event depending on $(\s(n))_{n=0}^N$, then
\be{5}
\P(A|S>N)\P(S>N) \leq \P(A) \leq \P(A|S>N) + \P(S \leq N),
\ee
so that  estimates for $\P(A)$ are obtained if estimates for $\P(A|S>N)$ and $\P(S>N)$ are available.

Similarly, to control that the system  performs {\em at most one} atypical updating {\em per time} up to time $N$, we define the random time
\be{3bis}
T := \min\left\{n \geq 1: \, \exists \, x \neq y \mbox{ such that } U_x(n) ,U_y(n) \not\in \left( \frac{e^{-2J+q}}{2\cosh(2J-q)}, 1- \frac{e^{-2J+q}}{2\cosh(2J-q)}\right) \right\}.
\ee

By definition $T\ge S$. We now establish  estimates for the random times $S$ and $T$ independently of the starting configurations.
From now on, when we need to indicate the initial condition $\s(0) = \s$, we write $\P_{\s}$ rather that $\P$ for the underlying probability.

We will adopt the following notation.
 For a given function $f:(0,+\infty) \ra (0,+\infty)$ we let $O(f(r))$ to be any function for which there is a constant $C>0$ satisfying $\frac{f(r)}{C} \leq O(f(r)) \leq C f(r)$ for $r \geq C$. Moreover, $a_r \sim b_r$ will stand for $\lim_{r \ra +\infty} \frac{a_r}{b_r} = 1$.

\bl{tempi}
There exist constants $C_i$ such that for each $a>0$ and $L$ sufficiently large we have
\be{9a}
\sup_{\s}\P_\s(S > L^{4k-2+a}) \leq C_1 e^{-O(L^a)}.
\ee
\be{11a}
\sup_{\s}\P_\s({S} \leq  L^{4k-2-a}) \leq C_2 L^{-a}
\ee
\be{8a}
\sup_{\s}\P_\s(T \leq L^{8k-4-a}) \leq C_3 L^{-a}
\ee
\be{8bisa}
\sup_{\s}\P_\s(T =S) \leq C_4 L^{-(4k-2)+2a}.
\ee
\el

{\em Proof.}
 To show \eqref{9a}, observe that $\{S>n\}$ means that up to time $n$ only {\em typical} updatings have been made. Since the probability that a given updating is typical is bounded above by $1- \frac{e^{-2J-q}}{2 \cosh(2J+q)} = 1-O(L^{-4k})$,
\[
\P_\s(S > L^{4k-2+a}) \leq \left(1-O(L^{-4k}) \right)^{L^2 \cdot L^{4k-2+a}} \leq C_1 e^{-O(L^a)},
\]
for some $C_1>0$, which establishes \eqref{9a}. To prove \eqref{11a}, observe that
\[
\begin{split}
\P_\s({S} \leq  L^{4k-2-a}) & = P\left( \exists x \in \L, \, n \leq L^{4k-2-a}: \, U_x(n) \not\in \left( \frac{e^{-2J+q}}{2\cosh(2J-q)}, 1- \frac{e^{-2J+q}}{2\cosh(2J-q)}\right) \right) \\
& \leq L^2 L^{4k-2-a}  \frac{e^{-2J+q}}{\cosh(2J-q)} = O(L^{-a}).
\end{split}
\]
The proof of \eqref{8a} is similar, the difference being that at least two atypical updatings need to occur:
\begin{multline*}
\P_\s(T \leq L^{8k-4-a}) \\ = P\left( \exists x,y \in \L, \, n \leq L^{8k-4-a}: \, U_x(n), U_y(n) \not\in \left( \frac{e^{-2J+q}}{2\cosh(2J-q)}, 1- \frac{e^{-2J+q}}{2\cosh(2J-q)}\right) \right) \\
 \leq L^4 L^{8k-4-a}  \left(\frac{e^{-2J+q}}{\cosh(2J-q)} \right)^2 = O(L^{-a}).
\end{multline*}
Finally, using \eqref{9a} and \eqref{8a},
\[
\begin{split}
\P_\s(T =S) & = \P_\s(T =S, \, S > L^{4k-2+a}) + \P_\s(T =S, \, T \leq L^{4k-2+a}) \\ &  \leq C_1 e^{-O(L^a)} + O(L^{-(4k-2)+2a}) = O(L^{-(4k-2)+2a}).
\end{split}
\]

\section{Mixing time and  tunneling time} \label{sec:mix}

In this section we prove Theorem \ref{mix} by giving estimates on the distribution of the hitting time
\[
T_{\bf 1} := \min\{n \geq 1: \s(n) = {\bf 1}\}.
\]
Since the dynamics described in the previous construction preserves the componentwise partial order on configurations, as noted in Remark \ref{monotonia}, we have
\be{6}
\sup_{\s \in \mathcal{S}}\P_{\s}(T_{\bf 1} > N) \leq \P_{\bf -1}(T_{\bf 1} > N).
\ee
 Thus, an upper bound on $\P_{\bf -1}(T_{\bf 1} > N)$ provides an upper bound for the mixing time. 
Indeed by using the coupling defined in Remark \ref{monotonia} we can define
the coupling time 
$$
\t_{couple}=\min\{n\ge 0:\; \s(n)=\s'(n)\}.$$
 The total variation distance between the evoluted measure at time $n$ and the stationary one, $d_{PCA}(n)$,  is related  to the coupling time  by the following
$$
d_{PCA}(n)\le \max_{\s,\s'}\P_{\s,\s'}(\t_{couple}>n)
$$
moreover, again due to the
monotonicity of the dynamics mentioned above, we have 
$$
\max_{\s,\s'}\P_{\s,\s'}(\t_{couple}>n)\le \P_{\bf -1}(T_{\bf 1} > n).
$$

So Theorem \ref{mix} immediately follows by the following:
 
\bt{th:main}
In the low temperature regime given in Definition \ref{lt-regime}, with $c > \frac12$ and $k-4c >4$,
\[
\lim_{L \ra +\infty} \P_{\bf -1}\left(T_{\bf 1} > L^{8k} \right) = 0.
\]
\et

The proof of this theorem is obtained in two steps and both are driven by the following idea.
We have three time scales given by three well separated order of magnitude of transition
probabilities. In the first scale the dynamics recurs in a very small subset of the state space ${\mathcal{S}}_1\subset{\mathcal{S}}$,
this recurrence can be described in terms of a suitable
1 dimensional random walk. 
On the second time scale the process jumps between different states in ${\mathcal{S}}_1$
and we can define a chain on this  restricted state space ${\mathcal{S}}_1$  and estimate
 its transition probabilities.
The third time scale is large enough with respect to the 
thermalisation of the random walk and thus can be ignored.

In the first step we show that,
due to the particular considered interaction, the
configurations with the same spin in each diagonal are stable
under the zero-temperature dynamics and  when the first atypical move
takes place, at time $S$, with large probability we have $S<T$ so that a single
discrepancy appears in a diagonal. The crucial remark is that starting with
such a configuration, the time 
$R$ needed to come back to diagonal configurations is typically much shorter than the waiting time for
the next atypical move,  so that
starting from ${\bf -1}$  the dynamics
can be studied in terms of a much simpler evolution moving in the space 
of diagonal configurations.

We need some notations. We denote by $\theta$ the {\em horizontal shift} on $\L$:
\[
\theta(i,j) = (i+1,j).
\]
By a common abuse of notation, we let $\theta$ act on configurations by $\theta \s_x := \s_{\theta(x)}$. For $m=0,1,L-1$, let $D_m$ denote the $m$-th NW-SE diagonal:
\[
D_m := \{ (i,j) \in \L: i+j = m\}
\]
(sums are, always, mod. $L$). Note that $D_{m+1} = \theta D_m$. The {\em diagonal configurations}, i.e. those that are constant on the diagonals, are denoted by:
\[
\mathcal{D} := \{ \s \in \mathcal{S} : \,x,y \in D_m \Rightarrow \s_x = \s_y \}.
\]
Assuming $\s(0) = \s \in \mathcal{D}$, it is immediately seen from the construction 
of the process given in Section \ref{2.1} that if only {\em typical} updatings occur up to time $N$, then $\s(n) = \theta^n \s$ for $n \leq N$. Thus, the evolution is trivial up to the stopping time $S$
and actually
\be{7}
S= \min\{n : \, \s(n) \neq \theta \s(n-1)\}.
\ee

Let $T$ be the time defined in \eqref{3bis}.
In the event $S < T$, which happens, as proven in Lemma \ref{tempi}  with high probability, $\s(S)$ is diagonal up to a single {\em discrepancy}, i.e. there is a unique $X \in \L$ such that $\s_X(S)$ is opposite to all other spins in the diagonal containing $X$, while $\s(S)$ is constant on all other diagonals.
Next Lemma shows that the site $X$ at which the first discrepancy appears is nearly uniformly distributed in $\L$.
\bl{lemma:1}

The conditional probability 
\[
 \P(X = x|  S < T)
 \]
is constant on both elements of the following partition of $\L$:
\[
\{x : \s_x = \s_x^l\}, \ \  \{x : \s_x = -\s_x^l\}
\]
and
\be{10}
\frac{e^{-4q}}{L^2} \leq \P(X = x|  S < T) \leq \frac{e^{4q}}{L^2} 
\ee
\el

The next step in our argument consists in studying the process from the time the first discrepancy appears to the next hitting time of $\mathcal{D}$, i.e. the time at which a diagonal configuration obtained. As we shall see, the time needed to go back to $\mathcal{D}$ is, with high probability, much shorter than the time needed for the next atypical updating to take place. 
 
For a rigorous analysis, under the condition $\{S < T\}$, we study the process $\{\s(S + n): n \geq 0\}$. By the strong Markov property, this is equivalent to study the process $\{\s(n): n \geq 0\}$  with an initial condition $\s(0) = \s$ which is diagonal, with a single discrepancy in $x \in D_m$.  Starting with such $\s$, besides typical and atypical updatings, neutral updatings arise. Indeed, the sites $x^u$ and $x^r$ can perform neutral updatings, having the left neighbor and the down neighbor of opposite sign. Suppose that no atypical updating occur. Then at time $1$ all diagonals are constant, except at most for the diagonal $D_{m+1}$. Here there are three possibilities:
\bi
\item[(i)]
both $\s^u_{x}$ and $\s_x^r$ update to $-1$: the discrepancy disappears, and $\s(1)$ is diagonal;
\item[(ii)]
both $\s^u_{x}$ and $\s_x^r$ update to $1$: the discrepancy has doubled, two neighboring sites in $D_{m+1}$ are $1$, while the rest of the diagonal is $-1$.
\item[(iii)]
in both other cases, the discrepancy has just shifted (up or right) to $D_{m+1}$. 
\ei
Under the condition of no atypical updatings, this argument can be repeated:  the discrepancy is shifted from a diagonal $D$ to $\theta D$, and its length can at most increase or decrease by one unit. The configuration goes back to $\mathcal{D}$ as soon as the discrepancies disappear or fill the whole diagonal. In order to keep fixed the diagonal containing the discrepancy, set
\[
\eta(n) := \theta^{-n} \s(n).
\]
If no atypical updating occur, $\eta$ remains constant except for the spins in $D_m$: here the number of spins equal to $1$ evolves as a random walk, that we show to be nearly symmetric. Standard estimates on random walks allow to estimate the probability the diagonal $D_m$ gets filled by ones before returning to all $-1$'s.

To make this argument precise, define the following stopping time:
\[
R := \min\{n>0 : \, \s(n) \in \mathcal{D} \}.
\]
Thus,  $R$ is the time the configuration has returned to $\mathcal{D}$.
\bl{lemma:2}
Assume the initial configuration $\s$ is diagonal with a single discrepancy at $x$,
i.e., if $x\in D_m$ then $\s_x=-\s_y$ for all $y\not=x$ in $D_m$, call ${\mathcal D}_x$ these configurations. Assume $2k-4c-3>0$. Then, for all $1<r<2k-1$
\be{12}
\begin{split}
\P_{\s}(R> L^r | {S} > L^{2k}) & \leq O(L^{-r+1}) \\ 
\P_{\s}(R> L^r ) & \leq O(L^{-r+1}).
\end{split}
\ee
\be{13}
\P_{\s}(\eta_x(R) = \s_x |{S} > L^{2k}) \sim \left\{ \begin{array}{ll} 4c \frac{\log L}{L} & \mbox{if } \s_x^u = \s_x \\ \frac{4c \log L}{L^{4c+1}} & \mbox{if }  \s_x^u = -\s_x \end{array} \right.
\ee
Moreover, let  $\eta^{D_m}$ be the configuration obtained from $\eta$ by flipping all spins in $D_m$. Then
\be{13bis}
\P_{\s}\left(\eta(R) = \eta^{D_m}\right) \sim \left\{ \begin{array}{ll} 4c \frac{\log L}{L} & \mbox{if } \s_x^u = \s_x \\ \frac{4c \log L}{L^{4c+1}} & \mbox{if }  \s_x^u = -\s_x \end{array} \right.
\ee
\el

Before continuing our argument, we comment on the meaning of these inequalities. Since by \eqref{11a} we know that the probability that an atypical updating occurs before time $L^{2k}$ is small, inequality \eqref{12} implies, in particular, that the configuration goes back to $\mathcal{D}$ in a time much shorter that ${S}$ (we are assuming $k$ large). Inequality \eqref{13} states that the probability the initial discrepancy at $x$ propagates to the whole diagonal is much higher if $\s_x$, $x \in D_m$, is equal to the spins in $D_{m+1}$. Most importantly, Lemma \ref{lemma:2} provides estimates on the transition for a starting diagonal configuration $\s \in \mathcal{D}$ to the {\em next} diagonal configuration hit after having left $\mathcal{D}$. This suggests to study an {\em effective process} obtained by observing $\eta(n)$ only at the times it enters $\mathcal{D}$.

Define the stopping times
\be{rem-times}
\begin{split}
R_0 & := 0 \\
S_1 & := \min\{m >0 : \, \s(m) \not\in \mathcal{D} \}=S \\
R_n & := \min\{m > S_n : \, \s(m) \in \mathcal{D} \}=S_n+R\circ\Theta_{S_n} \\
S_{n+1} & := \min\{m > R_n : \s(m) \not\in \mathcal{D} \}=R_n+S\circ\Theta_{R_n}
\end{split}
\ee
where $\Theta_t$ is the time shift operator acting on each trajectory of the Markov Chain
$\{\s(0),\s(1),....\}$ as a shift
$$
\Theta_t\{\s(0),\s(1),....\}=\{\s(t),\s(t+1),....\}.
$$
The following estimates follow from Lemma \ref{lemma:1} and Lemma \ref{lemma:2}.
\bc{cor:1}
The following estimates hold for all $n \geq 0$:
\be{16}
\P\left(S_{n+1} - R_n > L^{5k} \right) \leq e^{-L^k}.
\ee
\be{17}
\P\left(R_n - S_n > L^k \right) \leq \sup_{\s\in\cup_x{\mathcal D}_x} \P_\s\left(R>L^k\right)+\sup_\s \P_\s\left(S=T\right) \leq 
O(L^{-k+1}).
\ee
\ec

We now consider the Markov chain $(\eta(n))_{n \geq 0}$ at the times $R_n$ where the chain visits $\mathcal{D}$; more precisely we define
\be{reduced}
\xi(n) := \eta(R_n).
\ee
By the strong Markov property, $(\xi(n))_{n \geq 0}$ is a Markov chain in $\mathcal{D}$. Estimates on its transition probability are given in the following statement.
\bc{cor:2}
For all $\eta \in \mathcal{D}$ the following estimates hold.
\bi
\item[(a)]
If $\eta_x = - \eta_y$ for $x \in D_m$, $y \in D_{m+1}$ (we say $D_m$ is a {\em favorable} diagonal), then
\be{18}
\P\left(\xi(n+1) = \eta^{D_m} | \xi(n) = \eta\right) \geq  O\left(\frac{\log L}{L^2} \right).
\ee
Moreover, the above conditional probability is constant in $m$ on both elements of the partition of $\{0,1,\ldots,L-1\}$:
\[
\{m: x \in D_m, y \in D_{m-1} \Rightarrow \s_x = \s_y\}, \ \ \{m: x \in D_m, y \in D_{m-1} \Rightarrow \s_x = - \s_y\}.
\]
\item[(b)]
If $\eta_x =  \eta_y$ for $x \in D_m$, $y \in D_{m+1}$ ($D_m$ is an unfavorable diagonal), then
\be{19}
O\left(L^{-4c-2}\right) \leq \P\left(\xi(n+1) = \eta^{D_m} | \xi(n) = \eta\right) \leq  O\left(L^{-4c-1}\right).
\ee
\item[(c)]
\be{20}
\P\left(\xi(n+1) \not\in \{ \eta, \eta^{D_m}: \, m=0,\ldots,L-1\}| \xi(n) = \eta \right) \leq O(L^{-k+1})
\ee
\ei
\ec
\proof
Estimates \eqref{18} and \eqref{19} follow from \eqref{13bis} and the fact (see \eqref{10}) that a discrepancy is nearly uniformly distributed in $\L$ (Lemma \ref{lemma:1} ). Estimate \eqref{20} follows for the observation that  if $\xi(n+1) \not\in \{ \eta, \eta^{D_m}: \, m=0,\ldots,L-1\}$, then necessarily either two atypical updatings have occurred simultaneously between times $R_n$ and $S_{n+1}$, or an atypical updating have have occurred between times $S_{n+1}$ and $R_{n+1}$; the probability of this event has been estimated in (\ref{8bisa}) (used here with $a = k-1$) and (\ref{12}).
\qed

The process $\xi(n)$ defined in \eqref{reduced} starts at $\xi(0) = {\bf -1}$, and it can clearly identified with a process taking values in $\{-1,1\}^L$. Thus we write $\xi = (\xi_i)_{i=0}^{L-1}$, where $\xi_i$ is the spin in the diagonal $D_i$. By \eqref{19}, after  a waiting time of order at most $L^{4c+2}$, a one is created at some $i$. At this point there are two favorable diagonals: $D_i$ and $D_{i-1}$; all other diagonals are unfavorable. Thus,
in one time step, two transitions are equally likely: $\xi_i$ goes back to $-1$ or $\xi_{i-1}$ flips to $1$. By \eqref{18}, these transitions occur with probability $p \geq O\left(\frac{\log L}{L^2} \right)$. The probability that $\xi$ changes to some other configurations is, by \eqref{19} and \eqref{20}, not larger than $ O\left(L^{-k+1}\right) + O\left(L^{-4c-1}\right)$. In the case $\xi$ is back to ${\bf -1}$ the process starts afresh. Otherwise, there are two consecutive ones at $i-1,i$. The above argument can be iterated: in the next step two diagonals are favorable, $D_i$ and $D_{i-2}$, so $\xi_{i-2}$ and $\xi_i$ flips with the same probability $p$. Therefore, with overwhelming probability, the ones in $\xi(n)$ are consecutive, and their number evolves, up to events of small probability, as a symmetric $p$ random walk. This makes simple, for this effective process, to give estimates on the hitting time of ${\bf 1}$.

\bl{lemma:3}
Define $H^{(\xi)}_{\bf 1}$ the first time $\xi(n)$ visits $\{{\bf 1}\}$. Then, assuming $c > \frac12$ and $k-4c >4$,
\[
\P\left( H^{(\xi)}_{\bf 1} > L^{k+2} \right) \leq O(L^{-1}).
\]
\el
We are now ready to complete the proof of Theorem \ref{th:main}. Indeed using also Corollary \ref{cor:1},
\[
\P_{\bf -1}\left(T_{\bf 1} > L^{8k} \right) \leq \P\left(H^{(\xi)}_{\bf 1}>L^{k+2} \right) +
\P\left(R_{L^{k+2}}>L^{8k}\right)\leq
\]
\[ O(L^{-1})+
 \sum_{n \leq L^{k+2}} \P\left(R_{n} - R_{n-1} > L^{7k-2} \right) = O(L^{-1}),
\]
which is the desired result. 

\subsection{Proofs of the Lemmas} 
We are therefore left with the proof of Lemmas \ref{lemma:1},  \ref{lemma:2} and \ref{lemma:3}.

\noindent
{\em Proof of Lemma \ref{lemma:1}}

For the proof of \eqref{10}, recall that an atypical updating is made at $x$ at time $n$ if $U_n(x) \in I_x(\s(n-1))$, where
\[
I_x(\s) = \left\{ \begin{array}{ll} \left(0, \frac{e^{-2J + q \s_x}}{2 \cosh(-2J + q \s_x)}\right) & \mbox{if } \s_x^d = \s_x^l = -1 \\
\left(\frac{e^{2J + q \s_x}}{2 \cosh(2J + q \s_x)},1 \right) & \mbox{if } \s_x^d = \s_x^l = 1 \end{array} \right.
\]
We have:
$$
\P\left( X=x | T > S\right) =\frac{1}{\P(T>S)}\sum_n\P\left( X=x,\, S=n, \, T>n\right)
$$
and 
\[
\{X=x,\, S=n, \, T>n \} = \{S > n-1\} \cap \left\{ U_x(n)  \in I_x(\s(n-1)), \, U_y(n)  \not\in I_y(\s(n-1)) \mbox{ for } y \neq x \right\}.
\]
so that
$$
\P\left( X=x,\, S=n, \, T>n\right)=\P\left(  S>n-1\right)|I_x(\s(n-1))|\prod_{ y \neq x}\Big(1-|I_y(\s(n-1))|\Big)=
$$
$$ \P\left(  S>n\right)\frac{|I_x(\s(n-1))|}{1-|I_x(\s(n-1))|}=: \P\left(  S>n\right) f_x(n)
$$
We note that the function $f_x(n)$ as a function on $x$, is constant on the sets 
\[
M_+=\{x : \s_x = \s_x^l\}, \ \ M_-= \{x : \s_x = -\s_x^l\}, 
\]
so on these sets $\P\left( X=x | T > S\right) $ is constant, say $\P\left( X=x | T > S\right) =P_{M_{\pm}}$ .
Moreover since
$$
\min_\s |I_x(\s)|\ge e^{-4q}\max_\s |I_x(\s)|
$$
we have uniformly in $n$ 
$$
e^{-4q}<\frac{f_x(n)}{f_y(n)}<e^{4q}
$$ 
and so 
$$e^{-4q}<\frac{P_{M_+}}{P_{M_-}}<e^{4q},\qquad |M_+| P_{M_+}+|M_-|P_{M_+}=1$$
from which \eqref{10} easily follows.
\qed

\noindent
{\em Proof of Lemma \ref{lemma:2}}. 
We  prove \eqref{12} and \eqref{13}. The second inequality in \eqref{12} follows form the first, \eqref{11a} and the assumption $r <2k-1$, since
\[
\P_{\s}(R> L^r )\le\P_{\s}(R> L^r |S > L^{2k})  + \P_{\s}({S} \leq L^{2k} )
\]
 Note that, under $\P_{\s}(\cdot |{S} > L^{2k})$, the random numbers $\{U_x(n): \, x \in \L, \, n \leq L^{2k} \}$ are i.i.d. with uniform distribution on $\left( \frac{e^{-2J+q}}{2\cosh(2J-q)}, 1- \frac{e^{-2J+q}}{2\cosh(2J-q)}\right) $. The following probability describe the two possible neutral updatings; atypical updatings are forbidden by the conditioning.
\[
\begin{split}
\P(\eta_x(1) = 1 |{S} > L^{2k}) & = \P(\eta_{x^u}(1) = 1 |{S} > L^{2k})= \frac{e^{q \s_x^r}}{2 \cosh(q \s_x^r)} \\ & = \frac12 + \frac{c\s_x^r}{2} \frac{\log L}{L} + O\left(\left( \frac{\log L}{L} \right)^2 \right).
\end{split}
\]
Thus, denoting by $N(n)$ the number of spins equal to $1$ in the restriction to $D_m$ of $\eta(n)$, we have that 
\[
\begin{split}
p_+ := \P(N(1) = 2 |{S} > L^{2k}) & = \left(\P(\eta_x(1) = 1 |{S} > L^{2k})\right)^2 \\ & = \frac{1}{4} + \frac{c \s_x^r}{2} \frac{\log L}{L} + O\left(\left( \frac{\log L}{L} \right)^2 \right) \\
p_- := \P(N(1) = 0 |{S} > L^{2k}) & = \left(1- \P(\eta_x(1) = 1 |{S} > L^{2k})\right)^2 \\ & = \frac{1}{4} - \frac{c \s_x^r}{2} \frac{\log L}{L} + O\left(\left( \frac{\log L}{L} \right)^2 \right).
\end{split}
\]
This argument can now be repeated, since either the discrepancy for $\eta$ in $D_m$ has disappeared, or two neutral updatings are possible. This implies that, for $n \leq L^{2k}$ and $m>0$
\[
\begin{split}
p_+ & = \P(N(n) = m+1 |N(n-1) = m, \, {S} > L^{2k}) \\
p_- & = \P(N(n) = m-1 |N(n-1) = m, \, {S} > L^{2k})
\end{split}
\]
So, set $\tilde{R} := \min\{ n: N(n) \in \{0,L\}\}$. Note that $\tilde{R}\wedge L^{2k} = R\wedge L^{2k}$ on $\{T > S > L^{2k} \}$. Moreover, up to time $\tilde{R}\wedge L^{2k}$, $N(n)$ evolves as a $(p_+, p_-)$ one dimensional random walk. We recall that if $(\xi(n))_{n \geq 1}$ is a $(p_+, p_-)$ random walk with $\xi(0) = 1$, and denote by $H_{0L}, H_0, H_L$ the hitting times of, respectively, $\{0,L\}$, $\{0\}$ and $\{L\}$, then (see e.g. \cite{feller1}, XIV.2 and XIV.3, where the case $p_+ + p_- = 1$ is treated, but the same proof applies to $p_+ + p_- < 1$)
\be{mart1}
\P(H_L<H_0)=\frac{1-\frac{p_-}{p_+}}{1-\Big(\frac{p_-}{p_+}\Big)^L} \, \sim \, \left\{ \begin{array}{ll} \frac{4c \log L}{L} & \mbox{if } \s_x^u =\s_x \\
																\frac{4c \log L}{L^{4c+1}} & \mbox{if } \s_x^u = - \s_x \end{array} \right.
\ee
\be{mart2}
\E (H_{0L})=\frac{1}{p_+ - p_-}\Bigg[ L\frac{1-\frac{p_-}{p_+}}{1-\Big(\frac{p_-}{p_+}\Big)^L}-1 \Bigg]
\, \sim \, \left\{ \begin{array}{ll} 4L & \mbox{if } \s_x^r =1 \\
																\frac{ L}{c \log L} & \mbox{if } \s_x^r = - 1 \end{array} \right.
\ee
In particular, by Markov inequality, for every $r>1$
\be{15}
\P(H_{0L} > L^r) \leq O(L^{-r+1}).
\ee
From \eqref{mart1} and \eqref{15}, the desired estimate \eqref{12} and \eqref{13} follow. Finally, \eqref{13bis} follows from \eqref{9a} and\eqref{13}, using the assumption $2k-4c-3>0$.
\qed

\noindent
{\em Proof of Lemma \ref{lemma:3}}
Let 
\[
T^{\xi} := \min\{n: \xi(n) \not\in \{ \xi(n-1), \xi^{D_m}(n-1): \, m=0,\ldots,L-1\}\} .
\]
By \eqref{20}, 
\[
\P(T^{\xi}\leq L^{k-2}) \leq L^{k-2} O\left(L^{-k+1}\right) = O(L^{-1}).
\]
Similarly with what we did in previous Lemmas, we condition the Markov chain $\xi(n)$ to the event $\{T^{\xi} > L^{k-2}\}$. Under this conditioning, we are left with a Markov chain for which, up to time $L^{k-2}$, \eqref{18} and \eqref{19} hold, but transitions of the type in \eqref{20} are forbidden. Let
\[
S^{(\xi)}_1 := \min\{n: \, \xi(n) \neq {\bf -1} \}
\]
be the first time the process leaves the initial configuration, and
\[
\overline{S}^{(\xi)}_1 := \min\{n > S^{(\xi)}_1 : \xi(n) = \xi^i(n-1) \mbox{ for some } i \mbox{ such that } \xi_i(n-1) = \xi_{i+1}(n-1) \},
\]
where $\xi^i$ is the configuration obtained from $\xi$ by flipping $\xi_i$.
By \eqref{19}
\[
\P\left(S^{(\xi)}_1>L^{4c+3}\right) \leq  \left(1-  O\left(L^{-4c-2}\right)  \right)^{L^{4c+3}} \leq e^{-O(L)},
\]
and
\[
\P\left( \overline{S}^{(\xi)}_1 - S^{(\xi)}_1 \leq L^{2c} \right)  \leq L^{2c} L O\left(L^{-4c-1}\right)  = O\left(L^{-2c}\right).
\]
Conditioning to the event $\{ T^{\xi}>L^{k-2}, \, S^{(\xi)}_1 \leq L^{4c+3}, \overline{S}^{(\xi)}_1- S^{(\xi)}_1 > L^{2c} \}$ which, for $k-4c$ large enough, has probability at least $1-O\left(L^{-2c}\right) \geq 1-O(L^{-1})$, the number of spin equal to $1$ in $\xi(S^{(\xi)}_1 + n)$ evolves as a symmetric random walk, starting from $1$, and moving with probability $p \geq O\left(\frac{\log L}{L^2} \right)$ (see \eqref{18}). We now use identities analogous to \eqref{mart1} and \eqref{mart2} for the case  $p_+ = p_- = p$:
\be{hitL}
\P(H_L<H_0)= \frac{1}{L},
\ee
and
\[
\E (H_{0L}) = \frac{L-1}{2p} \leq O(L^3).
\]
It follows that
\[
\P\left( \xi(H_{0L}) = {\bf 1}, H_{0L}<\overline{S}^{(\xi)}_1 | T^{\xi}>L^{k-2}, \, S^{(\xi)}_1 \leq L^{4c+3}, \overline{S}^{(\xi)}_1- S^{(\xi)}_1 > L^{2c} \right) \geq O(L^{-1}),
\]
and
\[
\P(H_{0L} > C) \leq \frac{O(L^3)}{C}.
\]
Thus, introducing the stopping times, for $j \geq 1$ (note the analogy with \eqref{rem-times} in the previous step of the renormalization)
\[
\begin{split}
R^{(\xi)}_0 & := 0 \\ S^{(\xi)}_j & := \min\{n>R^{(\xi)}_{j-1}: \xi(n) \not\in \{{\bf -1}, {\bf 1}\} \} \\
R^{(\xi)}_{j} & := \min\{n>S^{(\xi)}_j : \xi(n) \in \{{\bf -1}, {\bf 1}\} \}
\end{split}
\]
we have, by \eqref{hitL},
\[
\P\left(\xi(R^{(\xi)}_j) = {\bf 1} | \xi(R^{(\xi)}_{j-1}) = {\bf -1} \right) \geq O(L^{-1}),
\]
and 
\[
\P\left(R^{(\xi)}_j - R^{(\xi)}_{j-1} > L^k \right) \leq \P\left(S^{(\xi)}_1 > L^{k-1}\right) + \P(H_{0L} > L^{k-1})  \leq O(L^{-k+4}),
\]
where we have used again the fact that $k-4c$ is sufficiently large.
Finally, for $k$ large enough,
\[
\begin{split}
\P\left(T^{(\xi)}_{\bf 1}>L^{k+2} \right) & \leq \P\left(R^{(\xi)}_{L^2} \leq L^{k} \right) + \sum_{j \leq L^2} \P\left(R^{(\xi)}_j - R^{(\xi)}_{j-1} > L^k \right) \\
& \leq \left(1- O(L^{-1}) \right)^{L^2} + L^2 O(L^{-k+4}) \leq O(L^{-1}).
\end{split}
\]
\qed

\end{document}